# Some properties of surfaces of finite *III*-type

HASSAN AL-ZOUBI
Department of Mathematics
Al-Zaytoonah University of Jordan
P.O. Box 130, Amman 11733
JORDAN
dr.hassanz@zuj.edu.jo

*Abstract:* - In this paper, we firstly investigate some relations regarding the first and the second Laplace operators corresponding to the third fundamental form *III* of a surface in the Euclidean space $E^3$. Besides, we introduce the finite Chen type surfaces of revolution with nonvanishing Gauss curvature with respect to the third fundamental form. We present a special case of this family of surfaces of revolution in $E^3$, namely, surfaces of revolution with *R* is constant, where *R* denotes the sum of the radii of the principal curvature of a surface.

## 1 Introduction

One of the most interesting and profound aspects of differential geometry is the idea of surfaces of finite type which was born by B-Y. Chen in the early 1970s and since then, it has become a source of interest for many researchers in this field. The reader can refer to [17] for more details. In the framework of this kind of study, the first-named author with S. Stamatakis has given in [28] a new generalization to this area of study by giving a similar definition of surfaces of finite type.

Let *x* be an isometric immersion of a surface *S* in the 3-dimensional Euclidean space $E^3$. We represent by $\Delta^J$ the Laplacian operator of *S* acting on the space of smooth functions $C^\infty(S)$. Then *S* is said to be of finite *J*-type, *J = I, II, III*, if the position vector *x* of *S* can be decomposed as a finite sum of eigenvectors of $\Delta^J$ of *S*, that is

$$x = x_0 + x_1 + x_2 + \ldots + x_k, \quad (1)$$

where
$$\Delta^J x_i = \lambda_i x_i, \; i = 1, \ldots, k,$$
$x_0$ is a fixed vector and $\lambda_1, \lambda_2, \ldots, \lambda_k$ are eigenvalues of the operator $\Delta^J$.

Surfaces of finite type in $E^3$ regarding the second fundamental form were investigated for some important classes of surfaces. More precisely, the class of ruled surfaces was studied in [7], while in [3], H. AL-Zoubi studied tubular surfaces in $E^3$. Other classes such as translation surfaces, Quadric surfaces, surfaces of revolution, helicoidal surfaces, cyclides of Dupin, and spiral surfaces, the classification of its finite *II*-type surfaces still unknown. According to the third fundamental form, ruled surfaces in [4], tubes in [5], and quadric surfaces [6] are the only classes were investigated in $E^3$.

This type of study can be also extended to any smooth map, not necessary for the position vector of the surface, for example, the Gauss map of a surface. Regarding this see [8, 9].

Another generalization to the above, one can study surfaces in $E^3$ whose position vector *x* satisfies the following condition

$$\Delta^J x = Ax, \quad J = I, II, III, \quad (2)$$

where $A \in \mathfrak{R}^{3\times 3}$.

Related to this, in [29] it was proved that the spheres and the catenoids are the only surfaces of revolution satisfying the above equation. Similarly, in [1] it was shown that helicoids and spheres are the only quadric surfaces in $E^3$ that satisfy (2). Next, in [2] condition (2) was studied for the class of translation surfaces. In fact, authored ascertained that Scherk's surface is the only translation surface in the Euclidean 3-space that satisfies (2), finally, in [24] the authors studied the class of translation surfaces in Sol$_3$ satisfying (2). Surfaces satisfying condition (2) are said to be of coordinate finite *J*-type.

Another interesting study is to find surfaces in $E^3$ whose Gauss map *N* satisfies the relation (2) that is

$$\Delta^J N = AN, \quad J = I, II, III,$$

For this problem, readers can be referred to [10, 11, 13, 18, 19, 20, 21, ].

Interesting research also one can follow the idea in [23,26] by defining the first and second Laplace operator using the definition of the fractional vector operators.

In order to achieve our goal, we briefly introduce a formula for $\Delta^{III} x$ and $\Delta^{III} N$ by using tensors calculations. Further, in the last section, we contribute to the solution of our main result.

## 2 Fundamentals

We consider a smooth surface $S$ in $E^3$ given by a patch $x = x(u^1, u^2)$ on a region U: $= (a, b) \times \mathcal{R}$ of $\mathcal{R}^2$ in which does not contain parabolic points. We denote by

$$I = g_{ij} du^i du^j, \quad II = b_{ij} du^i du^j, \quad III = e_{ij} du^i du^j$$

the three fundamental forms of $S$. For any two differentiable functions $f(u^1, u^2)$ and $g(u^1, u^2)$ on $S$, the first differential parameter of Laplace regarding the fundamental form $J$ is defined by [12]

$$\nabla^J(f,g) := d^{ij} f_{/i} g_{/j}, \quad (3)$$

where $f_{/i} := \dfrac{\partial f}{\partial u^i}$ and $(d^{ij})$ denotes the inverse tensor of $(g_{ij})$, $(b_{ij})$ and $(e_{ij})$ for $J = I, II$ and $III$ respectively.

We first prove the following relations:

$$\nabla^I(f, x) + \nabla^{II}(f, N) = 0, \quad (4)$$

$$\nabla^{II}(f, x) + \nabla^{III}(f, N) = 0. \quad (5)$$

For the proof of (5) we use (3) and the Weingarten equations

$$N_{/j} = -e_{jk} b^{km} x_{/m} = -b_{jk} g^{km} x_{/m}, \quad (6)$$

to obtain
$$\nabla^{II}(f, N) = b^{ij} f_{/i} N_{/j} = -b^{ij} f_{/i} b_{jk} g^{km} x_{/m}$$
$$= -g^{im} f_{/i} x_{/m} = -\nabla^I(f, x),$$

being (4). We have similarly

$$\nabla^{III}(f, N) = e^{ij} f_{/i} N_{/j} = -e^{ij} f_{/i} e_{jk} b^{km} x_{/m}$$
$$= -b^{im} f_{/i} x_{/m} = -\nabla^{II}(f, x),$$

which is (5).

The second Laplace operator according to the fundamental form $J = I, II, III$ of $S$ is defined by [10]

$$\Delta^J f := -d^{ij} \nabla_i^J f_{/j},$$

where $f$ is a sufficiently differentiable function, $\nabla_i^J$ is the covariant derivative in the $u^i$ direction with respect to the fundamental form $J$ [12]. For $J = III$ we have

$$\Delta^{III} f = -e^{ij} \nabla_i^{III} f_{/j}, \quad (7)$$

We now compute $\Delta^{III} x$ and $\Delta^{III} N$. From (7) and the equations [19, p.128]

$$\nabla_j^{III} x_{/i} = -b^{km} \nabla_m^I b_{ij} x_{/k} + b_{ij} N$$

we get

$$\Delta^{III} x = e^{ij} b^{km} \nabla_m^I b_{ij} x_{/k} - e^{ij} b_{ij} N. \quad (8)$$

Denote by
$$\Lambda_{ij}^k := \frac{1}{2} e^{km} (-e_{ij/m} + e_{im/j} + e_{jm/i}),$$

the Christoffel symbols of the second kind regarding the third fundamental form. We put

$$T_{ij}^k := \Gamma_{ij}^k - \Pi_{ij}^k.$$
$$\tilde{T}_{ij}^k := \Lambda_{ij}^k - \Pi_{ij}^k.$$

It is known that [19, p.22]

$$T_{ij}^k = -\frac{1}{2} b^{km} \nabla_m^I b_{ij}, \quad (9)$$

$$\tilde{T}_{ij}^k = -\frac{1}{2} b^{km} \nabla_m^{III} b_{ij} \quad (10)$$

and

$$T_{ij}^k + \tilde{T}_{ij}^k = 0. \quad (11)$$

Besides, using Ricci's Lemma
$$\nabla_j^{III} e^{ik} = 0$$

and the formula
$$R = \frac{2H}{K} = e^{ik} b_{ik}, \quad (12)$$

where $K$ is the Gauss curvature and $H$ is the mean curvature of $S$ respectively we have

$$R_{/m} = \nabla_m^{III}(e^{ik} b_{ik}) = e^{ik} \nabla_m^{III} b_{ik}. \quad (13)$$

From (9), (10), (11) and (13) we find

$$e^{ij} b^{km} \nabla^{I}_{m} b_{ij} = -2e^{ij} T_{ij}{}^{k} = 2e^{ij} \tilde{T}_{ij}{}^{k}$$
$$= -e^{ij} b^{km} \nabla^{III}_{m} b_{ij} = -b^{km} R_{/m}$$

and so

$$e^{ij} b^{km} \nabla^{I}_{m} b_{ij} \, x_{/k} = -b^{km} R_{/m} \, x_{/k} = -\nabla^{II}(R, x). \quad (14)$$

By combining (8), (12), and (14) we obtain [22]

$$\Delta^{III} x = -\nabla^{II}(R, x) - R N.$$

Finally, using (5) we arrive at

$$\Delta^{III} x = \nabla^{III}(R, N) - R N. \quad (15)$$

For the normal vector $N$ we have

(1.23) $\qquad \nabla^{III}_{k} N_{/i} = -e_{ik} N$

we have

$$\Delta^{III} N = -e^{ik} \nabla^{III}_{k} N_{/i} = e^{ik} e_{ik} N,$$

so that we conclude

$$\Delta^{III} N = 2 N.$$

From the last equation, it can be seen that the Gauss map of every surface $S$ in $E^3$ is of finite $III$-type 1, the corresponding eigenvalue is 2. Now we prove some relations according to the third fundamental form of f.

For any differentiable function $f(u^1, u^2)$ it can be easily shown that

$$\Delta^{III}(f x) = (\Delta^{III} f) x + f \Delta^{III} x - 2\nabla^{III}(f, x)$$
$$= (\Delta^{III} f) x + f \nabla^{III}(R, N) - fR N - 2\nabla^{III}(f, x)$$

Similarly

$$\Delta^{III}(f N) = (\Delta^{III} f) N + f \Delta^{III} N - 2\nabla^{III}(f, N)$$
$$= (\Delta^{III} f) N + 2f N + 2\nabla^{II}(f, x)$$

Denote by $W = -\langle x, N \rangle$ the support function of $S$, where $\langle , \rangle$ is the Euclidean inner product. Applying relation (7) for the function $W$, it can be easily verified that

$$\Delta^{III} W = -e^{ik} \nabla^{III}_{k} W_{/i} = e^{ik} \nabla^{III}_{k} \langle x, N \rangle_{/i}$$
$$= \langle e^{ik} \nabla^{III}_{k} x_{/i}, N \rangle + \langle x, e^{ik} \nabla^{III}_{k} N_{/i} \rangle$$

$$+ 2e^{ik} \langle x_{/i}, N_{/k} \rangle = -\langle \Delta^{III} x, N \rangle$$
$$- \langle x, \Delta^{III} N \rangle - 2e^{ik} \langle x_{/i}, b_{kt} g^{tr} x_{/r} \rangle$$
$$= -\langle N, \nabla^{III}(R, N) - RN \rangle$$
$$- \langle x, 2N \rangle - 2e^{ik} \langle x_{/i}, b_{kt} g^{tr} x_{/r} \rangle$$
$$= 2W - R. \quad (16)$$

We consider now the surface $S$ of finite $III$-type 1. Then we have $\Delta^{III} x = kx$, where $k$ is a constant eigenvalue.

From (15) we get $\nabla^{III}(R, N) - RN = kx$. Taking the inner product of both sides of this equation with $N$ we find $R = kW$, From the formula (16) we find that

$$\Delta^{III} W = (2 - k)W, \quad \Delta^{III} R = (2 - k)R.$$

Thus, we have proved the following:

**Theorem 1.** *Let S be a surface in $E^3$ of finite III-type 1 with corresponding eigenvalue l. Then the support function W and the sum of the principal radii of curvature R are of eigenfunctions of the Laplacian $\Delta^{III}$ with corresponding eigenvalue $2 - k$.*

Let now $S$ be a minimal surface. Then we have
$$R = \frac{2H}{K} = 0.$$

Thus from the equation (16) we get $\Delta^{III} W = 2W$. So we have

**Corollary 1**. *Let S be a minimal surface. Then the support function W is of an eigenfunction of $\Delta^{III}$ with corresponding eigenvalue equal 2.*

Let $S^*$ be a parallel surface of $S$ in (directed) distance $\mu = $ const. $\neq 0$, so that $1 - 2\mu H + \mu^2 K \neq 0$,. Then $S^*$ possesses the position vector $x^* = x + \mu N$.

Denoting by $K^*$ and $H^*$ the Gauss and mean curvature of $S^*$ respectively, we mention the following relations

$$K^* = \frac{K}{1 - 2\mu H + \mu^2 K},$$
$$H^* = \frac{H - \mu K}{1 - 2\mu H + \mu^2 K}.$$

Hence we get

$$R^* = \frac{2H^*}{K^*} = R - 2\mu \quad .(17)$$

On the other hand, the surfaces $S$, $S^*$ have common unit normal vector and spherical image. Thus $III = III^*$ and $\Delta^{III} = \Delta^{III*}$. We prove now the following theorem for later use.

**Theorem 2**. *Let $S$ be a minimal surface in $E^3$. Then $S^*$ is a parallel surface of $S$ if and only if the sum of the principal radii of curvature $R^*$ of $S^*$ is constant.*

*Proof.* Suppose that $S$ is a minimal surface in $E^3$, which is defined on a simply connected domain $D$ in the $(u^1, u^2)$-plane. Let

$$S^*: \boldsymbol{x}^* = \boldsymbol{x} + \mu \boldsymbol{N}, \mu \neq 0$$

be parallel surface of $S$. From (17) and taking into account $H = 0$, we find $R^* = -2\mu = \text{const}$. Hence the first part of the theorem is proved.

Conversely, let $R^* = \text{const}. \neq 0$. Then from Theorem (4.4) (see [29]), $S^*$ is of null $III$- type 2. Therefore from (1) there exist nonconstant eigenvectors $\boldsymbol{x}_1(u^1, u^2)$ and $\boldsymbol{x}_2(u^1, u^2)$ defined on the same domain $D$ such that

$$\boldsymbol{x}^* = \boldsymbol{x}_1 + \boldsymbol{x}_2, \quad (18)$$

where $\Delta^{III} \boldsymbol{x}_1 = \lambda_1 \boldsymbol{x}_1$, $\Delta^{III} \boldsymbol{x}_2 = \lambda_2 \boldsymbol{x}_2$, and here we have $\lambda_1 = 0$ because $S^*$ is of null $III$- type 2.

Once we have $\Delta^{III} \boldsymbol{x}^* = \Delta^{III} \boldsymbol{x}_1 + \Delta^{III} \boldsymbol{x}_2$, it then follows that

$$\Delta^{III} \boldsymbol{x}^* = \lambda_2 \boldsymbol{x}_2. \quad (19)$$

Besides, since $R^* = \text{const}. \neq 0$, we find

$$\Delta^{III} \boldsymbol{x}^* = -R^* \boldsymbol{N}. \quad (20)$$

Thus from (19) and (20), one finds

$$\lambda_2 \boldsymbol{x}_2 = -R^* \boldsymbol{N}$$

or $\boldsymbol{x}_2 = c \boldsymbol{N}$, where $c = -\dfrac{R^*}{\lambda_2}$, and then (18) becomes

$$\boldsymbol{x}^* = \boldsymbol{x}_1 + c\boldsymbol{N}. \quad (21)$$

The differential of the above equation is

$$d\boldsymbol{x}^* = d\boldsymbol{x}_1 + cd\boldsymbol{N}. \quad (22)$$

Taking the inner product of both sides of (22) with $\boldsymbol{N}$ yields

$$<d\boldsymbol{x}_1, \boldsymbol{N}> = 0. \quad (23)$$

Now we want to show that $\boldsymbol{x}_1(u^1, u^2)$ is a regular parametric representation of a surface in $E^3$. It is enough to prove that

$$\boldsymbol{x}_{1/1} \times \boldsymbol{x}_{1/2} \neq \boldsymbol{0}, \quad \forall (u^1, u^2) \in D,$$

where $\times$ is the Euclidean cross product. We have

$$\boldsymbol{x}_1 = \boldsymbol{x}^* - \mu \boldsymbol{N}^* \quad (24)$$

Using the Weingarten equations

$$\boldsymbol{N}_{/i} = -b_{ij} g^{jk} \boldsymbol{x}^*_{/k}.$$

and the equation (24), it follows that

$$\boldsymbol{x}_{1/1} \times \boldsymbol{x}_{1/2} = (\boldsymbol{x}^*_{/1} - \mu \boldsymbol{N}_{/1}) \times (\boldsymbol{x}^*_{/2} - \mu \boldsymbol{N}_{/2})$$

$$= (\boldsymbol{x}^*_{/1} \times \boldsymbol{x}^*_{/2}) - \mu(\boldsymbol{x}^*_{/1} \times \boldsymbol{N}_{/2}) + \mu(\boldsymbol{x}^*_{/2} \times \boldsymbol{N}_{/1})$$

$$+ \mu^2 (\boldsymbol{N}_{/1} \times \boldsymbol{N}_{/2})$$

$$= (1 - 2\mu H + \mu^2 K)(\boldsymbol{x}^*_{/1} \times \boldsymbol{x}^*_{/2}) \neq \boldsymbol{0},$$

$$\forall (u^1, u^2) \in D. \quad (25)$$

Hence, on account of (23) and (25), we conclude that $\boldsymbol{x}_1(u^1, u^2)$ is a regular parametric representation of a surface in $E^3$ with $\boldsymbol{N}$ its Gauss map.

Since $\Delta^{III} \boldsymbol{x}_1 = \boldsymbol{0}$. Consequently, from Theorem (3.1) (see [29]), $\boldsymbol{x}_1(u^1, u^2)$ is a minimal surface. Thus from (21), we obtain that $S^*$ is a parallel surface of a minimal. Now we mention and prove our main theorem.

**Theorem 3.** *The only surfaces of revolution in $E^3$ of which the sum of the radii of the principal curvature $R$ is constant are*
• *parts of spheres which are of finite III-type 1,*
• *catenoids which are of finite null III-type 1, and*
• *the parallel surfaces to the catenoids, which are of finite null III-type 2.*

## 3  Proof of Theorem 3

Let $C$ be a smooth curve lies on the $xz$-plane parametrized by

$$\boldsymbol{x}(u) = (f(u), 0, g(u)), u \in J, (J \subset \mathcal{R}),$$

where $f$, $g$ are smooth functions and $f$ is a positive function. When $C$ is revolved about the $z$-axis, the resulting point set $S$ is called the surface of revolution generated by the curve $C$. In this case, the $z$-axis is called the axis of revolution of $S$, and $C$ is called the profile curve of $S$. On the other hand, a subgroup of the rotation group which fixes the vector $(0, 0, 1)$ is generated by

$$\begin{pmatrix} cosv & -sinv & 0 \\ sinv & cosv & 0 \\ 0 & 0 & 1 \end{pmatrix}.$$

Then the position vector of $S$ is given by see ([14, 24])

$$x(u, v) = (f(u) \cos v, f(u) \sin v, g(u)), \quad (26)$$

$$u \in J, v \in [0, 2\pi).$$

Without loss of generality, we may assume that $C$ has the arc-length parametrization, i.e., it satisfies

$$(f')^2 + (g')^2 = 1, \quad (27)$$

where $' := \frac{d}{du}$. Furthermore if $f' g' = 0$, then $f = $ const. or $g = $ const. and $S$ would be a circular cylinder or part of a plane, respectively. A case that has been excluded since $S$ would consist only of parabolic points.

The partial derivatives of (26) are

$$x_u = (f'(u) \cos v, f'(u) \sin v, g'(u)),$$

and

$$x_v = (-f(u) \sin v, f(u) \cos v, 0).$$

The components $g_{ij}$ of the first fundamental form in (local) coordinates are the following

$$g_{11} = <x_u, x_u> = 1, \ g_{12} = <x_u, x_v> = 0,$$
$$g_{22} = <x_v, x_v> = f^2.$$

Denoting by $\kappa$ the curvature of the curve $C$ and $r_1$, $r_2$ the principal radii of curvature of $S$, we have

$$r_1 = \frac{1}{\kappa}, \qquad r_2 = \frac{f}{g'}.$$

The Gauss curvature and the mean curvature of $S$ are respectively

$$K = \frac{1}{r_1 r_2} = \frac{\kappa g'}{f} = -\frac{f''}{f},$$

$$2H = \frac{1}{r_1} + \frac{1}{r_2} = \kappa + \frac{g'}{f}.$$

The Gauss map $N$ of $S$ is computed as follows

$$N(u,v) = (-g'\cos v, -g'\sin v, -f'). \quad (28)$$

Now, by using the natural frame $\{N_u, N_v\}$ of $S$ defined by

$$N_u = (-g'' \cos v, -g'' \sin v, f'')$$

and

$$N_v = (g' \sin v, -g' \cos v, 0)$$

the components $e_{ij}$ of the third fundamental form in (local) coordinates are the following

$$e_{11} = <N_u, N_u> = (g'')^2 + (f'')^2,$$

$$e_{12} = <N_u, N_v> = 0, \ e_{22} = <N_v, N_v> = (g')^2.$$

The Beltrami operator $\Delta^{III}$ in terms of local coordinates $(u, v)$ of $S$ can be expressed as follows

$$\Delta^{III} = -\frac{1}{\kappa^2}\frac{\partial^2}{\partial u^2} + \frac{g'\kappa' - \kappa g''}{g'\kappa^3}\frac{\partial}{\partial u} - \frac{1}{g'^2}\frac{\partial^2}{\partial v^2}. \quad (29)$$

On account of (27) we put

$$f' = \cos \varphi, \ g' = \sin \varphi,$$

where $\varphi = \varphi(u)$. Then $\kappa = \varphi'$ and the parametric representation (28) of the unit vector $N$ of $S$ becomes

$$N(u, v) = \{-\sin\varphi \cos v, -\sin\varphi \sin v, \cos\varphi \}. \quad (30)$$

Also relation (29) takes the following form

$$\Delta^{III} = -\frac{1}{\varphi'^2}\frac{\partial^2}{\partial u^2} + \left(\frac{\varphi''}{\varphi'^3} - \frac{\cos\varphi}{\varphi'\sin\varphi}\right)\frac{\partial}{\partial u} - \frac{1}{\sin^2\varphi}\frac{\partial^2}{\partial v^2}. \quad (31)$$

For the sum of the principal radii of curvature $R = r_1 + r_2 = \frac{2H}{K}$, one finds

$$R = \frac{f}{\sin \varphi} + \frac{1}{\varphi'}. \quad (32)$$

Taking the derivative of (32) we find

$$R' = -\frac{\varphi''}{\varphi'^2} - \frac{f\varphi'\cos\varphi}{\sin^2\varphi} + \frac{\cos\varphi}{\sin\varphi}. \quad (33)$$

Let $(x_1, x_2, x_3)$ be the coordinate functions of (26). By virtue of (31), we obtain

$$\Delta^{III} x_1 = \Delta^{III}(f \cos v) = \left( \frac{\varphi'' \cos \varphi}{\varphi'^3} - \frac{1}{\varphi' \sin \varphi} + \frac{2 \sin \varphi}{\varphi'} + \frac{f}{\sin^2 \varphi} \right) \cos v \tag{34}$$

$$\Delta^{III} x_2 = \Delta^{III}(f \sin v) = \left( \frac{\varphi'' \cos \varphi}{\varphi'^3} - \frac{1}{\varphi' \sin \varphi} + \frac{2 \sin \varphi}{\varphi'} + \frac{f}{\sin^2 \varphi} \right) \sin v \tag{35}$$

$$\Delta^{III} x_3 = \Delta^{III} g = -\frac{2 \cos \varphi}{\varphi'} + \frac{\varphi'' \sin \varphi}{\varphi'^3}. \tag{36}$$

From (32) and (33), equations (34), (35) and (36) become respectively

$$\Delta^{III} x_1 = \left( R \sin \varphi - \frac{R' \cos \varphi}{\varphi'} \right) \cos v, \tag{37}$$

$$\Delta^{III} x_2 = \left( R \sin \varphi - \frac{R' \cos \varphi}{\varphi'} \right) \sin v, \tag{38}$$

$$\Delta^{III} x_3 = -\frac{R' \sin \varphi}{\varphi'} - R \cos \varphi. \tag{39}$$

We obtain the following two cases:

*Case I.* $R \equiv 0$. Thus $H \equiv 0$. Consequently $S$, being a minimal surface of revolution, is a catenoid.

*Case II.* $R = \text{const.} \neq 0$. From (37), (38), and (39) we obtain

$$\left. \begin{array}{l} \Delta^{III} x_1 = R \sin \varphi \cos v \\ \Delta^{III} x_2 = R \sin \varphi \sin v \\ \Delta^{III} x_3 = -R \cos \varphi \end{array} \right\} \tag{40}$$

Let $(N_1, N_2, N_3)$ be the coordinate functions of $N$. From (29), relations (40) can be written

$$\Delta^{III} x_1 = -R N_1, \quad \Delta^{III} x_2 = -R N_2, \quad \Delta^{III} x_3 = -R N_3,$$

and hence
$$\Delta^{III} \mathbf{x} = -R \mathbf{N}. \tag{41}$$

In view of (7) and (41) we have

$$(\Delta^{III})^m \mathbf{x} = -(2^{m-1}) R \mathbf{N}. \tag{42}$$

Now, if $S$ is of finite type $k$, then there exist real numbers $c_1, c_2, \ldots, c_k$ such that

$$(\Delta^{III})^k \mathbf{x} + c_1 (\Delta^{III})^{k-1} \mathbf{x} + \ldots + c_k \mathbf{x} = \mathbf{0}. \tag{43}$$

From (41) and (42) relation (43) becomes

$$-2^{k-1} R \mathbf{N} - 2^{k-2} c_1 R \mathbf{N} - \ldots - c_{k-1} R \mathbf{N} + c_k \mathbf{r} = \mathbf{0},$$
or
$$c \mathbf{N} + c_k \mathbf{x} = \mathbf{0}, \tag{44}$$

where $c = -R(2^{k-1} + 2^{k-2} c_1 + \ldots + c_{k-1}) = \text{const.}$.

Now, if $c_k \neq 0$, then from (44) we have $\mathbf{x} = -\frac{c}{c_k} \mathbf{N}$, and hence we get $|\mathbf{x}| = \left| \frac{c}{c_k} \right|$ and so $S$ is a sphere. On account of Theorem (3.3) (see [23]), $S$ is of finite $III$-type 1. If $c_k = 0$, then $S$ is of null type $k$. Since $R = \text{const.}$, thus according to Theorem (4.4) (see [29]) and Theorem (2), $S$ is of null $III$-type 2 which is a parallel surface of a minimal.